# SHAPES OF AUSLANDER-REITEN TRIANGLES

EDSON RIBEIRO ALVARES, SÔNIA MARIA FERNANDES, AND HERNÁN GIRALDO

*Dedicated to Professor Hector Merklen on the occasion of his 80th birthday*

ABSTRACT. Our main theorem classifies the Auslander-Reiten triangles according to properties of the morphisms involved. As a consequence, we are able to compute the mapping cone of an irreducible morphism. We finish by showing a technique for constructing the connecting component of the derived category of any tilted algebra. In particular we obtain a technique for constructing the derived category of any tilted algebra of finite representation type.

## 1. INTRODUCTION

Let $\Lambda$ be a finite dimensional non semisimple algebra over a field $k$. The notion of irreducible morphism and almost split sequence (Auslander-Reiten sequence) in the category of finitely generated right modules over $\Lambda$, denoted by $\mathrm{mod}\,\Lambda$, was introduced by M. Auslander [4] and M. Auslander and I. Reiten [5] and it plays an important role in the study of the category $\mathrm{mod}\,\Lambda$. So, the Auslander-Reiten theory allows us to better understand the module category of an Artin algebra, and the Auslander-Reiten quiver is a visualization of this module category. From this quiver, one may know all the isomorphism classes of indecomposable modules in $\mathrm{mod}\,\Lambda$, if $\Lambda$ is of finite representation type, on may also know all the homomorphism between the indecomposable modules. The derived category of the module category of a finite-dimensional algebra can be very helpful for understanding the module category. For example, there are various invariants associated with an algebra $\Lambda$ that are preserved under Morita equivalence but not only by Morita equivalence, but also by derived equivalence of the module category.

The theory of almost split sequences was successfully extended to the context of derived category of a finite global dimensional algebra, where the notion of almost split sequence gave rise to the notion of almost split triangles, that will be called Auslander-Reiten triangles here. However, not all triangulated category of a finite dimensional algebra has Auslander-Reiten triangles. In ( [15], [13]) Happel proved that the bounded derived category of a finite dimensional algebra has Auslander-Reiten triangles if and only if the global dimension of the algebra is finite. In this paper, he also proved that the first two morphisms in an Auslander-Reiten triangle are irreducible morphisms and described the Auslander-Reiten triangle of the derived category of a finite dimensional hereditary algebra. He was very successful in giving conditions to an Auslander-Reiten sequence to become an Auslander-Reiten triangle. There are many reasons to investigate the Auslander-Reiten triangles and irreducible morphisms.

The structure of the Auslander-Reiten triangles in the particular case of full subcategory of compact objects was studied by Beligiannis in [8]. An investigation was also carried out on the connection between Auslander-Reiten triangles and the existence of Serre functors and the representability of homological functors. In [21], Reiten-Van den Bergh demonstrated that the existence of a right Serre functor is equivalent to the existence of right Auslander-Reiten triangles, in triangulated $\mathrm{Hom}$-finite Krull-Schmidt $k$-categories. Gao has examined this in a particular class of algebras [11]. He studied the conditions for the Gorenstein derived categories of finite-dimensional Gorenstein algebras to admit Auslander-Reiten triangles. Krause has also been investigating the existence of Auslander-Reiten triangle in [19].

They have been studied by Wheeler, aiming to describe the Auslander-Reiten quiver of the derived category [26]. In this direction, Scherotzke shows that, if the Auslander-Reiten quiver of the derived category has a finite component whose elements belong to the bounded homotopy category of finitely generated projectives, then the algebra itself is simple and the component is $A_1$ (see [24]). Another result was obtained by Scherotzke in this paper. Diveris, Purin and Webb have extended a previous result of this paper in [10].

Roggenkamp tried to give a unified approach to Auslander-Reiten triangles in the case of a finite dimensional algebra over a field, as well as in the case of a classical order over a complete Dedekind domain (see [23]). Liu investigated a unified approach to the study of almost split sequences (in abelian categories) and Auslander-Reiten triangles (in triangulated categories) through the concept of an Auslander-Reiten sequence in a Krull-Schmidt category ([20]). Bobiński has described all Auslander-Reiten triangles in the category of perfect complexes of gentle algebras explicitly [9]. This class of algebras was introduced by I. Assem and A. Skowroński [2] whose indecomposable complexes of bounded derived category have been classified by V. Bekkert and H. A. Merklen in [7]. The proof relies on D. Happel's embedding functor from the bounded derived category to the stable module category of the repetitive algebra, and on the known description of Auslander-Reiten triangles in this stable module category. The same strategy was carried out in another investigation of Happel, in which he proved that, if $f$ is an irreducible morphism in the derived category of a finite dimensional algebra, then its embedding in the stable module category is an irreducible morphism [16].

In the classical Auslander-Reiten theory for module categories over Artin algebras, subcategories (usually not abelian subcategories) of $\mathrm{mod}\,\Lambda$ having Auslander-Reiten sequences were studied by M. Auslander and S. O. Smalø (see [3]). Jørgensen investigated in [18] the analogue for triangulated categories.

Recent studies aim to extend the concept of irreducible morphisms in $\mathrm{mod}\,\Lambda$ to irreducible morphisms in the derived category of $\mathrm{mod}\,\Lambda$. Giraldo and Merklen [12] studied irreducible morphisms in the categories $\mathcal{C}(\mathcal{A})$ and $\mathrm{D}^-(\mathrm{mod}\,\Lambda)$, where $\mathcal{C}(\mathcal{A})$ is the category of complexes over an abelian Krull-Schmidt category $\mathcal{A}$ and $\mathrm{D}^-(\mathrm{mod}\,\Lambda)$ is the derived category of the bounded above complexes of finite generated right $\Lambda$-modules. A study on Auslander-Reiten sequences for complexes with applications to the bounded derived category over an artinian algebra was given by Bautista-Salorio in [6]. The irreducible morphism endings in perfect complexes, and the irreducible morphisms between bounded complexes of projectives in case of a self injective algebra were throughly analysed.

This paper aimed to investigate the shape of Auslander-Reiten triangles. We have shown that an Auslander-Reiten triangle that begins with an irreducible morphism of type smonic, has to finish in an irreducible morphism of type sepic. An Auslander-Reiten triangle that begins with an irreducible morphism of type sepic, has to finish with an irreducible morphism of type sirreducible and if an Auslander-Reiten triangle begins with an irreducible morphism of type sirreducible, then it has to finish either as smonic or sirreducible. We go one step further with this characterization. We classify the shape of the mapping cone in each case and show a special shape for each one of them. We shall give here examples that demonstrate that these three conditions (smonic, sepic or sirreducible) are not enough to be irreducible.

The paper is organized as follows. In Section 2, we give some background material from [12], including properties of irreducible morphisms and Auslander-Reiten theory for triangulated categories. We keep the notation and conventions for categories of complexes over an additive category. For the convenience of the reader, we present a simple lemma about irreducible morphisms and give examples showing that, in general, smonic, sepic and sirreducible morphisms are not necessarily irreducible. In Section 3, we show the standard form of the smonic, sepic and sirreducible morphism that was presented in [12] and present the shape of the mapping cone of the smonic, sepic and sirreducible morphism in the homotopy category. We present a lemma showing that a morphism which is isomorphic for a smonic (sepic, sirreducible) morphism is a smonic (sepic, sirreducible, respectively) morphism. In Section 4, we classify the Auslander-Reiten triangle that begins with a smonic, sepic and sirreducible morphism. We also show that the complexes of each one of this Auslander-Reiten triangles assume a distinguished shape, which was thoroughly described in Theorem (2). In Section 5, we have shown that the mapping cone is indecomposable in a triangulated category. We also presented some particular results of the mapping cone of a irreducible morphism that begins in a module category, besides some particular conditions to the non existence of irreducible morphisms. Using the shapes developed in Section 4 we provides a technique, in Section 6, for constructing the connecting component of the derived category of any tilted algebra.

## 2. Basic facts

Let us first establish some notations. Given a quiver $\Gamma = (\Gamma_0, \Gamma_1)$, where $\Gamma_0$ is the set of vertices and $\Gamma_1$ is the set of arrows, $\alpha: i \to j$ denotes an arrow from the vertex $i$ to the vertex $j$. The notations $P_i$, $I_i$ and $S_i$ are used, respectively, for the indecomposable projective, injective and simple associated to the vertex $i \in \Gamma_0$.

$\mathcal{A}$ denotes an additive category (which in most of the applications will be an abelian, Krull-Schmidt, locally bounded k-category, where $k$ is a field) and $\mathrm{Hom}_{\mathcal{A}}(X, Y)$ denotes the set of all morphisms from $X$ to $Y$, where $X$ and $Y$ are objects in $\mathcal{A}$. We denote by $\mathcal{P}$ a full Krull-Schmidt subcategory of $\mathcal{A}$; and by $\mathcal{I}$ a (left, right or two-sided) ideal of $\mathcal{A}$. If $\mathcal{A}$ is the category of modules of an algebra, the set of all morphisms with image contained in the radical of the target is an ideal of $\mathcal{A}$ (see [12]).

We denote by $k$ a field and $\Lambda$ a finite dimensional $k$-algebra non semisimple with Jacobson radical $\mathbf{r}$. By $\mathrm{mod}\,\Lambda$ and $\mathcal{P}(\Lambda)$, we denote the category of finitely generated right $\Lambda$-modules and the full subcategory of finitely generated projective right $\Lambda$-modules, respectively.

We denote by $\mathcal{C}(\mathcal{A})$ the category of complexes $X = (X^i, d^i)_{i \in \mathbb{Z}}$ (resp. $\mathcal{C}(\mathcal{P})$ or $\mathcal{C}(\mathcal{P}(\Lambda))$) with cells in $\mathcal{A}$ (resp. in $\mathcal{P}$ or $\mathcal{P}(\Lambda)$); and we denote by $\mathcal{C}_{\mathcal{I}}(\mathcal{A})$ (resp. $\mathcal{C}_{\mathcal{I}}(\mathcal{P})$ or $\mathcal{C}_{\mathcal{I}}(\mathcal{P}(\Lambda))$) the category of those complexes whose differential morphisms are in $\mathcal{I}$ (resp. in $\mathcal{I} \cap \mathcal{P}$). We say that a complex $X$ is a minimal projective complex if $X \in \mathcal{C}_{\mathcal{I}}(\mathcal{P})$. By $\mathcal{K}_{\mathcal{I}}(\mathcal{A})$, we denote the homotopy category of $\mathcal{C}_{\mathcal{I}}(\mathcal{A})$; and by $\mathcal{K}_{\mathcal{I}}(\mathcal{P})$, the full subcategory of $\mathcal{K}_{\mathcal{I}}(\mathcal{A})$ whose complexes have cells in $\mathcal{P}$.

If $f: X \to Y$, $f = (f^n)_{n \in \mathbb{Z}}$, is a morphism of $\mathcal{C}(\mathcal{A})$, then $[f]$ denotes its homotopy class, the corresponding morphism of $\mathcal{K}(\mathcal{A})$. We usually say, that $f^n$ refers to the components of $f$ and also, by abuse of language, of $[f]$.

Throughout this article, we use the following usual notation for a morphism $f$ from a complex $X$ to a complex $Y$, $f = (f^n)_{n \in J}$, where $J = ]-\infty, m]$, $J = [n, m]$, $J = [n, \infty[$, or $J = \mathbb{Z}$ and $n < m$ are numbers in $\mathbb{Z}$. By $\mathcal{C}^J(\mathcal{P})$ (resp. $\mathcal{C}_{\mathcal{I}}^J(\mathcal{P})$) we denote the full subcategory of $\mathcal{C}(\mathcal{P})$ (resp. $\mathcal{C}_{\mathcal{I}}(\mathcal{P})$) determined by the complexes whose cells are zero outside $J$.

If $\mathcal{A}$ is an additive category and $f: X \to Y$ is a morphism in $\mathcal{A}$, $f$ is said to be a split monomorphism (resp. split epimorphism) if there is a morphism $h: Y \to X$ in $\mathcal{A}$ such that $hf = 1_X$ (resp. $fh = 1_Y$). When either one of these conditions hold, $f$ is said to be split. A morphism $f$ in $\mathcal{A}$ is said to be irreducible if it is not split and if for any factorization $f = hg$, either $h$ is a split epimorphism or $g$ is a split monomorphism.

When $\mathcal{A}$ is an abelian Krull-Schmidt category, $f: X \to Y$ is called radical if for any split monomorphism $g: X_1 \to X$ in $\mathrm{Hom}(X_1, X)$ and any split epi $h: Y \to Y_1$ in $\mathrm{Hom}(Y, Y_1)$, the composition $hfg$ is not an isomorphism. We denote by $\mathrm{rad}(X, Y)$ the set of all radical morphisms in $\mathrm{Hom}(X, Y)$. Note that a morphism $f: P \to Q$ between two projective finitely generated modules $P$ and $Q$, over an Artin algebra, is radical if and only if $Im f \subset rQ$.

If $\mathcal{P}(\Lambda)$ is the full subcategory of the category of right projective $\Lambda$-modules and $\mathcal{I}$ is the ideal consisting of the morphisms whose images are inside the radical of the target, $\mathcal{C}_{\mathcal{I}}(\mathcal{P}(\Lambda))$ is then the set of the so-called minimal projective complexes, ie, those complexes whose families of differential morphisms $d^n: P^n \to P^{n+1}$ are such that, for all $n$, the image of $d^n$ is contained in $\mathrm{rad}\, P^{n+1}$.

A complex $X$ is said to be bounded below if $X^i = 0$ for all but finitely many $i < 0$ and it is called bounded above if $X^i = 0$ for all but finitely many $i > 0$. Also, $X$ is bounded if it is bounded below and above. Finally, $X$ is said to have bounded cohomology if $H^i(X) = 0$ for all but finitely many $i \in \mathbb{Z}$. We denote by $\mathcal{C}^-(\mathcal{A})$ (resp. $\mathcal{C}^+(\mathcal{A})$, $\mathcal{C}^b(\mathcal{A})$) the full subcategory of complexes whose objects are the complexes bounded above (resp. bounded below, resp. bounded above and below). We denote by $\mathcal{C}^{-,b}(\mathcal{A})$ (resp. $\mathcal{C}^{+,b}(\mathcal{A})$) the full subcategory of complexes bounded above with bounded cohomology (resp. bounded below with bounded cohomology). We denote by $\mathcal{K}^{-,b}(\mathcal{A})$, $\mathcal{K}^{+,b}(\mathcal{A})$, and $\mathcal{K}^b(\mathcal{A})$ the homotopy categories of the categories of complexes introduced above.

If $\mathcal{A}$ has enough projective objects, then the localization $\mathrm{D}^-(\mathcal{A})$ of $\mathcal{K}^-(\mathcal{A})$ exists and is equivalent to the full subcategory $\mathcal{K}^-(\mathcal{P})$ of bounded above complexes of projective objects in $\mathcal{A}$. So, the derived category $\mathrm{D}^-(\mathrm{mod}\,\Lambda)$ is equivalent to $\mathcal{K}^-(\mathcal{P}(\Lambda))$ (Theorem 10.4.8 [25]). We describe the shape of the Auslander-Reiten triangles in $K^-(\mathcal{P}(\Lambda))$, so in the $\mathrm{D}^-(\mathrm{mod}\,\Lambda)$ because they are triangle equivalents. Hence, in order to use our theorems to characterize the Auslander-Reiten triangles of $\mathrm{D}^-(\mathrm{mod}\,\Lambda)$, we assume that each object of this category may be defined, up to isomorphism, by a complex $X$ which is a minimal projective complex in $\mathcal{K}^-(\mathcal{P}(\Lambda))$ and that the irreducible morphisms of $\mathcal{K}^-(\mathcal{P}(\Lambda))$ are given by the irreducible morphism of $\mathcal{C}_\mathcal{I}^-(\mathcal{P}(\Lambda))$. It is well known that any complex $X \in \mathcal{C}(\mathcal{P}(\Lambda))$ is isomorphic in $\mathcal{K}(\mathcal{P}(\Lambda))$ to a minimal projective complex.

Another important remark is that any morphism $f \in \mathcal{C}_\mathcal{I}^-(\mathcal{P})$ is irreducible if and only if $[f]$ is irreducible in $\mathcal{K}_\mathcal{I}^-(\mathcal{P})$ (Theorem 6 in [12]). When dealing with irreducible morphisms of complexes, we will always assume that at least one of the complexes is indecomposable. Let $f : X \to Y$ be a morphism of complexes. Then, clearly, if $f$ is a split monomorphism (resp. a split epimorphism), then $f^i$ is a split monomorphism (resp. epimorphism) for all $i \in J$. As observed in [12] (Example 2), this does not hold true in $\mathcal{K}(\mathcal{A})$.

2.1. **The basic definitions and results of irreducible morphisms.** We will give a definition below of three types of morphisms of complexes (see [12]) which will be used throughout this work.

**Definition 1.** *A morphism of complexes $f = (f^i)_{i \in J} : X \to Y$ in $\mathcal{C}(\mathcal{A})$ is called **smonic** (resp. **sepic**) if all its components, $f^i$, are split monomorphisms (resp. split epimorphisms) and is called **sirreducible** if there is exactly one index $\iota_0$ such that $f^{\iota_0}$ is irreducible in $\mathcal{A}$ and $f^i$ is a split epimorphism for $i < \iota_0$ and a split monomorphism for $i > \iota_0$.*

From now on, we refer to these three types as the *standard forms*. We will use the Proposition 3 and Corollary 2 of [12]. For the remainder, or to facilitate the exposition to the reader of this article, we have ennounced the result below.

**Proposition 1.** *(Proposition 3 - [12] ) Let $f = f^J : X_J \to Y_J$ be an irreducible morphism of complexes in $\mathcal{C}_\mathcal{I}^J(\mathcal{P})$, with $J$ an interval of $\mathbb{Z}$. One and only one of the next conditions hold:*

*(a) $f$ is a smonic morphism;*
*(b) $f$ is a sepic morphism;*
*(c) $f$ is a sirreducible morphism.*

The standard form for the first and the second kind of morphism was presented in Proposition 1 in [12]. For the convenience of the reader, we present in Proposition (2) the standard form for the three cases above.

Directly from the classification of the irreducible morphisms in $\mathcal{C}_\mathcal{I}^J(\mathcal{P})$ (Proposition 3 [12]) and in $\mathcal{K}_\mathcal{I}^J(\mathcal{P})$, we can state the following result.

**Lemma 1.** *Let $f : X \to Y$ be a morphism of complexes in $\mathcal{C}_\mathcal{I}^J(\mathcal{P})$.*

*(a) If $f$ is an irreducible morphism such that $X^i = 0$, for all $i < -t$, $X^{-t} \neq 0$ and $Y^{-k} \neq 0$ for some $k > t$, then $f$ is smonic.*
*(b) If $f$ is an irreducible morphism such that $X^j = 0$ for all $j > 0$, $Y^j = 0$ for all $j < -t$, $X^{-(t+1)} \neq 0$ and $Y^j \neq 0$ for some $j > 0$, then $f$ is sirreducible.*
*(c) If $f$ is an irreducible morphism such that $Y^j = 0$ for all $j < -t$ $(t > 0)$ and $j > 0$, $X^{-(t+1)} \neq 0$ and $X^1 \neq 0$, then $f$ is sepic.*
*(d) If $X^j = 0$ for all $j < -t$ (where $t > 0$), $X^1 \neq 0$, $Y^{-(t+1)} \neq 0$ and $Y^1 = 0$, then $f$ is not irreducible.*

2.2. **Being smonic, sepic or sirreducible is not enough to be an irreducible morphism.** We would like to point out using these examples that a smonic, sepic or sirreducible morphism of complexes, are not necessarily an irreducible morphism. In the next example, we will use the following notation for a $\Lambda$-module $M$ in the category $\mathcal{K}^b(\mathcal{P}(\Lambda))$ of an algebra $\Lambda$. If $0 \to P_n \to \cdots \to P_1 \to P_0 \to M \to 0$ is the minimal projective resolution of $M$, then the $kth$-shift of $M$ will be represented by $[P_n - \cdots - P_0][k]$ in $\mathcal{K}^b(\mathcal{P}(\Lambda))$, where $P_0$ is in degree $-k$ of the complex.

**Example 1.** *Let $\Lambda$ be the finite dimensional algebra given by the quiver*

$$1 \xleftarrow{\phantom{xx}} 2 \xleftarrow{\phantom{xx}} 3.$$

*The Auslander-Reiten quiver of $\mathcal{K}^b(\mathcal{P}(\Lambda))$ is the following:*

$$\cdots P_1[0] \cdots\cdots\cdots (P_1 - P_2 - P_3)[-1] \cdots\cdots\cdots P_3[0] \cdots\cdots\cdots P_2[1] \cdots$$

(diagram with morphisms $\epsilon, \gamma, \alpha, \zeta, \delta, \beta$ between objects $(P_2 - P_3)[-1], (P_1 - P_2)[0], (P_2 - P_3)[0], P_3[-1], P_2[0], P_1[1], (P_1 - P_2 - P_3)[0]$)

*Note that $\beta\alpha$ is smonic, $\delta\gamma$ is sepic and $\zeta\epsilon$ is sirreducible. However, it is easy to see that they are not irreducible morphisms.*

2.3. **The Auslander-Reiten theory for triangulated categories.** First, observe that the homotopy category and the derived category are triangulated categories ([25]), where the automorphism is the shift functor $[1]$. The distinguished triangles are (up to isomorphism) given by $X \xrightarrow{f} Y \xrightarrow{\iota_f} C_f \xrightarrow{p_f} X[1]$ for any morphism $f$ (where $C_f$ is the cone of $f$). This triangle will be called the standard triangle associated with $f$ in $\mathcal{K}(\mathcal{P})$.

Let $\mathcal{T}$ be a triangulated category with translation functor $T$. A distinguished triangle $X \xrightarrow{u} Y \xrightarrow{v} Z \xrightarrow{w} TX$ is called an Auslander-Reiten triangle if the following conditions are satisfied:

(a) The objects $X, Z$ are indecomposable.
(b) The morphism $w$ is non-zero.
(c) If $f : W \to Z$ is not a split epimorphism, then there exists $f' : W \to Y$ such that $vf' = f$.

According to Happel, we have that $(1) + (2) + (3)$ is equivalent to the condition $(1) + (2) + (3')$, where $(3')$: If $f : X \to W$ is not a section, then there exists $f' : Y \to W$ such that $f'u = f$. We say that the Auslander-Reiten triangle starts in $X$, has middle term $Y$ and ends in $Z$ and we refer to $w$ as the connecting homomorphism of an Auslander-Reiten triangle. Also note that an Auslander-Reiten triangle is uniquely determined up to isomorphism, and $u$ and $v$ are irreducible morphisms (Proposition 4.3 [14]).

The conditions for the existence of Auslander-Reiten triangles in a triangulated category have been determined in [21]. It is demonstrated that a triangulated category admits Auslander-Reiten triangles (that is, for every indecomposable element $X$, there is an Auslander-Reiten triangle that ends in $X$ and one that starts in $X$) if and only if the category has a Serre functor. In [13], it was proved that $D^b(\text{mod } \Lambda)$ has Auslander-Reiten triangles if the global dimension of $\Lambda$ is finite. After that, the following generalization was given in [15].

**Theorem**(Happel) Let $Z \in \mathcal{K}^{-,b}(\mathcal{P}(\Lambda))$ be indecomposable. Then, there is an Auslander-Reiten triangle ending in $Z$, if and only if $Z \in \mathcal{K}^b(\mathcal{P}(\Lambda))$.

It is well known that a morphism $f : X \to Y$ between indecomposable complexes in $\mathcal{C}^{-,b}(\mathcal{P}(\Lambda))$ is irreducible in $\mathcal{C}^{-,b}(\mathcal{P}(\Lambda))$ if and only if $f$ is irreducible in $\mathcal{K}^{-,b}(\mathcal{P}(\Lambda))$. We can therefore choose, for an irreducible morphism in $\mathcal{K}^{-,b}(\mathcal{P}(\Lambda))$, an irreducible morphism $\mathcal{C}^{-,b}(\mathcal{P}(\Lambda))$ that represents this morphism. We can restrict to the study of certain subcategories of finite complexes as shown in [6].

## 3. Standard Forms

For the convenience of the reader, we present below the result Proposition 1 [12]). In addition, we present the standard form of the sirreducible morphism to facilitate the exhibition of the next results. The proof of this third formulation follows in the same way of the other two cases.

### 3.1. The standard forms of the smonic, sepic and sirreducible morphisms.

**Proposition 2. [Standard Forms]**
Let $(X, d^i)_{i \in \mathbb{Z}}$ and $(Y^i, \partial^i)_{i \in \mathbb{Z}}$ be complexes and $f : X \to Y$ be a morphism of complexes in $\mathcal{C}^J(\mathcal{A})$ for an abelian category $\mathcal{A}$

$$
\begin{array}{ccccccccc}
\cdots & \longrightarrow & X^{i-1} & \xrightarrow{d^{i-1}} & X^i & \xrightarrow{d^i} & X^{i+1} & \longrightarrow & \cdots \\
& & \downarrow f^{i-1} & & \downarrow f^i & & \downarrow f^{i+1} & & \\
\cdots & \longrightarrow & Y^{i-1} & \xrightarrow{\partial^{i-1}} & Y^i & \xrightarrow{\partial^i} & Y^{i+1} & \longrightarrow & \cdots.
\end{array}
$$

(a) If $f$ is a smonic morphism, (up to isomorphism) we can and will assume that $Y^i = X^i \oplus Y'^i$, $f^i = \begin{pmatrix} 1 \\ 0 \end{pmatrix}$, and $\partial^i = \begin{pmatrix} d^i & a^i \\ 0 & e^i \end{pmatrix}$.

(b) If $f$ is a sepic morphism, (up to isomorphism) we will write $X^i = Y^i \oplus X'^i$, $f^i = (\,1\ 0\,)$, and $d^i = \begin{pmatrix} \partial^i & 0 \\ b^i & e^i \end{pmatrix}$.

(c) If $f$ has a irreducible component $f^{i_0}$ ($i_0 \in \mathbb{Z}$) and $f^n$ is a split epimorphism for $n < i_0$ and $f^n$ is an split monomorphism for $n > i_0$, then we can and will assume that
   (1) for all $i < i_0$, $X^i = Y^i \oplus X'^i$, $f^i = (\,1\ 0\,)$,
   (2) for $i < i_0 - 1$, we assume $d^i = \begin{pmatrix} \partial^i & 0 \\ b^i & e^i \end{pmatrix}$,
   (3) for $i_0$, $d^{i_0-1} = (\, c^{i_0-1}\ e^{i_0-1}\,) : Y^{i_0-1} \oplus X'^{i_0-1} \to X^{i_0}$ and $\partial^{i_0} = \begin{pmatrix} \ell^{i_0} \\ \epsilon^{i_0} \end{pmatrix} : Y^{i_0} \to X^{i_0+1} \oplus Y'^{(i_0+1)}$,
   (4) for $i > i_0 + 1$, we assume that $Y^i = X^i \oplus Y'^i$, $f^i = \begin{pmatrix} 1 \\ 0 \end{pmatrix}$, $\partial^i = \begin{pmatrix} d^i & a^i \\ 0 & \epsilon^i \end{pmatrix}$.

*Proof.* It is an easy verification. □

### 3.2. The standard triangles that begin in a sirreducible morphism.
In this section, we discuss the mapping cone of the sirreducible morphism in $\mathcal{K}^J(\mathcal{P}(\Lambda))$.

**Proposition 3.** Let $f : X \to Y$ be a morphism of complexes in $\mathcal{C}^J(\mathcal{P}(\Lambda))$ and $X \xrightarrow{f} Y \xrightarrow{\iota_f} C_f \xrightarrow{p_f} X[1]$ the standard triangle $\mathcal{K}^J(\mathcal{P}(\Lambda))$. If $f$ is sirreducible and for some $i$, $f^i : X^i \to Y^i$ is the unique irreducible morphism of $f$ and $f$ is in the standard form, then the standard triangle is isomorphic to the triangle $X \xrightarrow{f} Y \xrightarrow{g} Z \xrightarrow{w} X[1]$ in the homotopic category $\mathcal{K}^J(\mathcal{P}(\Lambda))$, where

$$
Z^j = \begin{cases} X'^{j+1} & \text{if } j < i-1, \\ X^i, & \text{if } j = i-1, \\ Y^i, & \text{if } j = i, \\ Y'^i, & \text{if } j > i \end{cases} \text{ with } d_Z^j = \begin{cases} -e^{j+1} & \text{if } j \leq i-2, \\ f^i, & \text{if } j = i-1, \\ \epsilon^j, & \text{if } j \geq i \end{cases}
$$

and $g$ and $w$ are the morphisms of complexes with

$$
g^j = \begin{cases} b^j & \text{if } j < i-1, \\ c^{i-1}, & \text{if } j = i-1, \\ 1, & \text{if } j = i, \\ (\,0\ 1\,), & \text{if } j > i \end{cases} \text{ and } w^j = \begin{cases} \begin{pmatrix} 0 \\ 1 \end{pmatrix} & \text{if } j < i-1, \\ 1, & \text{if } j = i-1, \\ -\ell^i, & \text{if } j = i, \\ -a^j, & \text{if } j > i \end{cases}
$$

*Proof.* We have that if $j < i-2$, $d_Z^j d_Z^{j-1} = e^{j+1}e^j = 0$ because $f$ is a morphism of complex. If $j = i-2$, $d_Z^{i-2} d_Z^{i-3} = e^{i-1}e^{i-2} = 0$ because $f$ is a morphism of complex. In the same way, if $j = i-1$, $d_Z^j d_Z^{j-1} = -f^i e^{i-1} = 0$ because $f^i d^{i-1} = 0$. If $j > i$, we have $d_Z^j d_Z^{j-1} = \epsilon^j \epsilon^{j-1} = 0$ because $\partial^j \partial^{j-1} = 0$. So, $Z$ is a complex. Now, the equalities $f^i (\, c^{i-1}\ e^{i-1}\,) = \partial^{i-1} (\,1\ 0\,)$, $(\, c^{i-1}\ e^{i-1}\,) \begin{pmatrix} \partial^{i-2} & 0 \\ b^{i-2} & e^{i-2} \end{pmatrix} = 0$, $\begin{pmatrix} \partial^j & 0 \\ b^j & e^j \end{pmatrix} \begin{pmatrix} \partial^{j-1} & 0 \\ b^{j-1} & e^{j-1} \end{pmatrix} = 0$ for $j \leq i-2$ give us the following commutative diagram

$$\cdots \xrightarrow{\partial^{i-3}} Y^{i-2} \xrightarrow{\partial^{i-2}} Y^{i-1} \xrightarrow{\partial^{i-1}} Y^i \xrightarrow{\partial^i} X^{i+1} \oplus Y'^{(i+1)} \xrightarrow{\partial^{i+1}} X^{i+2} \oplus Y'^{(i+2)} \to \cdots$$
$$\downarrow b^{i-2} \quad \downarrow c^{i-1} \quad \downarrow 1 \quad \downarrow (0\ 1) \quad \downarrow (0\ 1)$$
$$\cdots \xrightarrow{-e^{i-2}} X'^{(i-1)} \xrightarrow{-e^{i-1}} X^i \xrightarrow{f^i} Y^i \xrightarrow{\epsilon^i} Y'^{(i+1)} \xrightarrow{\epsilon^{i+1}} Y'^{(i+2)} \longrightarrow \cdots$$

and so $g$ is a morphism of complex. Now, the commutativity of the following diagram shows us that $w$ is a morphism of complexes

$$\cdots \longrightarrow X'^{(i-2)} \xrightarrow{-e^{(i-2)}} X'^{(i-1)} \xrightarrow{-e^{i-1}} X^i \xrightarrow{f^i} Y^i \xrightarrow{\epsilon^i} Y'^{(i+1)} \to \cdots$$
$$\downarrow \binom{0}{1} \quad \downarrow \binom{0}{1} \quad \downarrow 1 \quad \downarrow -\ell^i \quad \downarrow -a^{i+1}$$
$$\cdots \to Y^{i-2} \oplus X'^{(i-2)} \xrightarrow{d^{i-2}} Y^{i-1} \oplus X'^{(i-1)} \xrightarrow{d^{i-1}} X^i \xrightarrow{-d^i} X^{(i+1)} \xrightarrow{d^{i+1}} X^{(i+2)} \to \cdots$$

where $d^j = \begin{pmatrix} \partial^j & 0 \\ b^j & e^j \end{pmatrix}$ for $j < i-1$ and $d^{i-1} = \begin{pmatrix} c^{i-1} & \epsilon^{i-1} \end{pmatrix}$. To prove the isomorphism of triangles, we will define $h : C_f \to Z$, such that the following diagram is commutative:

$$\begin{array}{ccc} Y & \xrightarrow{\iota_f} C_f & \xrightarrow{p_f} X[1] \\ \downarrow 1 & \downarrow h & \downarrow 1 \\ Y & \xrightarrow{g} Z \xrightarrow{w} X[1] \end{array}.$$

If we define $h : C_f \to Z$ as $h^j = \begin{cases} (0\ 1\ b^j) & \text{if } j < i-1, \\ (1\ c^j), & \text{if } j = i-1, \\ (0\ 1), & \text{if } j = i, \\ (0\ 0\ 1), & \text{if } j > i, \end{cases}$ then we have that the diagram

$$Y^{i-2} \oplus X'^{(i-2)} \oplus Y^{i-1} \to Y^{i-1} \oplus X'^{(i-1)} \oplus Y^{i-2} \xrightarrow{d^{i-2}_{C_f}} X^i \oplus Y^{i-1} \xrightarrow{d^{i-1}_{C_f}} X^{i+1} \oplus Y^i \xrightarrow{d^i_{C_f}} X^{i+2} \oplus X^{i+1} \oplus Y'^{(i+1)}$$
$$\downarrow \qquad \downarrow (0\ 1\ b^{i-2}) \quad \downarrow (1\ c^{i-1}) \quad \downarrow (0\ 1) \quad \downarrow (0\ 0\ 1)$$
$$X'^{(i-2)} \xrightarrow{-e^{i-2}} X'^{(i-1)} \xrightarrow{-e^{i-1}} X^i \xrightarrow{f^i} Y^i \xrightarrow{\epsilon^i} Y'^{(i+1)}$$

with $d^{i-2}_{C_f} = \begin{pmatrix} -c^{i-1} & -\epsilon^{i-1} & 0 \\ 1 & 0 & \partial^{i-2} \end{pmatrix}, d^{i-1}_{C_f} = \begin{pmatrix} -d^i & 0 \\ f^i & \partial^{i-1} \end{pmatrix}, d^i_{C_f} = \begin{pmatrix} -d^{i+1} & 0 \\ 1 & \ell^i \\ 0 & \epsilon^i \end{pmatrix}$ is commutative. So, $h$ is a morphism of complexes. It is easy to see that $h\iota_f = g$. Now, we shall prove that $wh = p_f \circ 1$ in the homotopic category $\mathcal{K}^J(\mathcal{P}(\Lambda))$.

If we define $s^j = \begin{cases} \begin{pmatrix} 0 & 0 & 1 \\ 0 & 0 & 0 \\ 0 & 0 & 0 \end{pmatrix} & \text{if } j < i-1, \\ \begin{pmatrix} 0 & 0 \\ 0 & 1 \end{pmatrix}, & \text{if } j = i-1, \\ (0\ 0), & \text{if } j = i, \\ (0\ 1\ 0), & \text{if } j > i \end{cases}$ then we have that $p_f = wh$ in $\mathcal{K}^J(\mathcal{P}(\Lambda))$. Now, to prove that $h : C_f \to Z$ is an isomorphism in $\mathcal{K}^J(\mathcal{P}(\Lambda))$, let $\eta : Z \to C_f$, such that

$$s^j = \begin{cases} \begin{pmatrix} 0 \\ 1 \\ 0 \end{pmatrix} & \text{if } j < i-1, \\ \begin{pmatrix} 1 \\ 0 \end{pmatrix}, & \text{if } j = i-1, \\ \begin{pmatrix} -\ell^i \\ 1 \end{pmatrix}, & \text{if } j = i, \\ \begin{pmatrix} -a^j \\ 0 \\ 1 \end{pmatrix}, & \text{if } j > i. \end{cases}$$

It is easy to see that $\eta$ is a morphism of complexes and from $h^i \eta^i = (0\ 1) \begin{pmatrix} -\ell^i \\ 1 \end{pmatrix} = 1$ and $h^j \eta^j = 1$ for all $j \neq i, i-1$. So, $h\eta = 1$ in $\mathcal{C}^J(\mathcal{P}(\Lambda))$. Now, we prove that $\eta h = 1$ in $\mathcal{K}^J(\mathcal{P}(\Lambda))$. If we define $v^j : C_f^j \to C_f^{j-1}$ as

$$v^j = \begin{cases} \begin{pmatrix} 0 & 0 & 1 \\ 0 & 0 & 0 \\ 0 & 0 & 0 \end{pmatrix} & \text{if } j < i-1, \\ \begin{pmatrix} 0 & 0 & 0 \\ 1 & 0 & 0 \\ 0 & 0 & 0 \end{pmatrix}, & \text{if } j = i-1, \\ \begin{pmatrix} 0 & 0 \\ 0 & 0 \end{pmatrix}, & \text{if } j = i, \\ \begin{pmatrix} 0 & 1 & 0 \\ 0 & 0 & 0 \end{pmatrix}, & \text{if } j = i+1, \\ \begin{pmatrix} 0 & 1 & 0 \\ 0 & 0 & 0 \\ 0 & 0 & 0 \end{pmatrix}, & \text{if } j > i+1 \end{cases}$$

Then, it is easy to see that $\eta h = 1$ in $\mathcal{K}^J(\mathcal{P}(\Lambda))$. $\square$

### 3.3. The standard triangles that begin in a smonic morphism. 
In this section, we discuss the mapping cone of the smonic morphism in $\mathcal{K}^J(\mathcal{P}(\Lambda))$.

**Proposition 4.** *Let $f : X \to Y$ be a morphism of complexes in $\mathcal{C}^J(\mathcal{P}(\Lambda))$ and $X \xrightarrow{f} Y \xrightarrow{\iota_f} C_f \xrightarrow{p_f} X[1]$ the standard triangle associated with $f$ in $\mathcal{K}^J(\mathcal{P}(\Lambda))$. If $f$ is smonic morphism in the standard form, then the standard* *triangle is isomorphic to the triangle* $X \xrightarrow{f} Y \xrightarrow{g} W \xrightarrow{-a} X[1]$ *in the homotopic category* $\mathcal{K}^J(\mathcal{P}(\Lambda))$ *with* $W^i = Y'^i$, $d_W^i = e^i : Y'^i \to Y'^{i+1}$, $g^i = (\,0\ 1\,) : X^i \oplus Y'^i \to Y'^i$ *and* $(-a)^i = -a^i$.

*Proof.* The proof is similar to that already given. □

3.4. **The standard triangles that begin in a sepic morphism.** In this section, we discuss the mapping cone of the sepic morphism in $\mathcal{K}^J(\mathcal{P}(\Lambda))$.

**Proposition 5.** *Let* $f : X \to Y$ *be a morphism of complexes in* $\mathcal{C}^J(\mathcal{P}(\Lambda))$ *and* $X \xrightarrow{f} Y \xrightarrow{\iota_f} C_f \xrightarrow{p_f} X[1]$ *the standard triangle in* $\mathcal{K}^J(\mathcal{P}(\Lambda))$. *If* $f$ *is sepic morphism in the standard form, then the standard triangle is isomorphic to the triangle* $X \xrightarrow{f} Y \xrightarrow{g} Z \xrightarrow{\binom{0}{1}} X[1]$ *in the homotopic category* $\mathcal{K}^J(\mathcal{P}(\Lambda))$ *where* $Z$ *is the complex* $Z^i = X'^{i+1}, d_Z^i = -e^{i+1} : X'^{i+1} \to X'^{i+2}$ *and* $g^i = b^i : Y^i \to X'^{i+1}$.

*Proof.* The proof is similar to that already given. □

The next result allows us to show that the properties of the triangles in Propositions (3), (4) and (5) describe the general shape of the Auslander-Reiten triangles in the derived category.

3.5. **The sirreducible, smonic and sepic property is invariant by isomorphism.** From Corollary 3 and Lemma 4 of [12], we know that if $f \in \mathcal{C}_\mathcal{I}^J(\mathcal{P}(\Lambda))$ and $[f]$ is the homotopy class in $\mathcal{K}_\mathcal{I}^J(\mathcal{P}(\Lambda))$, then $f$ is split monomorphism (resp. split epimorphism) if and only if $[f]$ is split monomorphism (resp. split epimorphism). Now, we would like to prove that isomorphic irreducible morphisms in $\mathcal{K}_\mathcal{I}^J(\mathcal{P}(\Lambda))$ are in the same class, that is, if one of them is smonic (resp. sepic, sirreducible), then the other is smonic (resp. sepic, sirreducible).

**Lemma 2.** *Let* $f : X \to Y, g : X' \to Y' \in \mathcal{K}_\mathcal{I}^J(\mathcal{P}(\Lambda))$ *morphisms, such that* $f$ *is isomorphic to* $g$ *in* $\mathcal{K}_\mathcal{I}^J(\mathcal{P}(\Lambda))$.

*(a) If* $f$ *is smonic, then* $g$ *is smonic.*
*(b) If* $f$ *is sepic, then* $g$ *is sepic.*
*(c) Assume that* $f$ *and* $g$ *are irreducible. If* $f$ *is sirreducible, then* $g$ *is sirreducible.*

## 4. Shape of the Auslander-Reiten triangles

In general, it is difficult to compute the Auslander-Reiten translation. The shape of Auslander-Reiten triangles of an hereditary finite dimensional $k$-algebra was presented in [14]. Here, we show the general shape of Auslander-Reiten triangles in a derived category. There are three general cases. Essentially, in each case, we characterize the mapping cone of the Auslander-Reiten triangles and the kind of irreducible morphism in this triangle.

4.1. **A classification of the irreducible morphisms in the Auslander-Reiten triangles.**

**Theorem 1.** *Let* $X \xrightarrow{u} Y \xrightarrow{v} Z \xrightarrow{w} X[1]$ *be an Auslander-Reiten triangle in* $\mathcal{K}_\mathcal{I}^J(\mathcal{P}(\Lambda))$.
*(a) If* $u$ *is smonic, then* $v$ *is sepic.*
*(b) If* $u$ *is sepic, then* $v$ *is sirreducible.*
*(c) If* $u$ *is sirreducible, then* $v$ *is smonic or sirreducible.*

*Proof.* (a) If $u$ is smonic morphism in $\mathcal{C}_\mathcal{I}^J(\mathcal{P}(\Lambda))$, $u$ is isomorphic to a morphism $f$ in the standard form of the Proposition (2) (a). So, by Proposition (4), there exists a triangle $X \xrightarrow{f} Y \xrightarrow{g} W \xrightarrow{-a} X[1]$ in the homotopic category $\mathcal{K}_\mathcal{I}^J(\mathcal{P}(\Lambda))$ with $W^i = Y'^i$, $d_W^i = e^i : Y'^i \to Y'^{i+1}$, $g^i = (0,1) : X^i \oplus Y'^i \to Y'^i$ and $(-a)^i = -a^i$. So, this is a triangle in $\mathcal{K}_\mathcal{I}^J(\mathcal{P}(\Lambda))$ and so the triangle $X \xrightarrow{u} Y \xrightarrow{v} Z \xrightarrow{w} X[1]$ is isomorphic to the triangle $X \xrightarrow{f} Y \xrightarrow{g} W \xrightarrow{-a} X[1]$ in $\mathcal{K}_\mathcal{I}^J(\mathcal{P}(\Lambda))$ with $g$ sepic. By Lemma (2), we have that $v$ is sepic.

(b) Now, we suppose that $u$ is a sepic morphism in $\mathcal{C}_\mathcal{I}^J(\mathcal{P}(\Lambda))$, so $u$ is isomorphic to a morphism $f$ in the standard form of the Proposition (2) (b). So, by Proposition (5), if $f$ is sepic morphism in the standard form, then the standard triangle is isomorphic to the triangle $X \xrightarrow{f} Y \xrightarrow{g} Z \xrightarrow{\binom{0}{1}} X[1]$ in the homotopic category $\mathcal{K}_\mathcal{I}^J(\mathcal{P}(\Lambda))$

where $Z$ is the complex $Z^i = X'^{i+1}, d_Z^i = -e^{i+1} : X'^{i+1} \to X'^{i+2}, g^i = b^i : Y^i \to X'^{i+1}$. We have that $Im \, b^i \subseteq \text{rad } X'^{i+1}$. We know that $g$ is irreducible. We can say that $g$ is smonic, sepic or sirreducible. So, if each $b^i$ is split epimorphism, then $Im \, b^i = X'^i = \text{rad } X'^i$. So $X'^i = 0$ for all $i$. Then $Z = 0$ and so $f$ is an isomorphism. But we know that $f$ is not an isomorphism. Now if $b^i$ is split monomorphism for all $i$, then $Im \, \binom{0}{1} = Y^i$ and this is included in $\text{rad } Y^i$. So, $Y^i$ is zero for all $i$. However, $b^i$ is a split monomorphism, then $X'^{i+1} = 0$. Thus, $Z$ is zero, and we have that $f$ is an isomorphism, a contradiction.

(c) The proof is quite similar to the proof previously given. $\square$

### 4.2. The shape of the Auslander-Reiten triangles.

**Theorem 2.** *Let* $(*) \ X \xrightarrow{u} Y \xrightarrow{v} Z \xrightarrow{w} X[1]$ *be an Auslander-Reiten triangle in* $\mathcal{K}_\mathcal{I}^J(\mathcal{P}(\Lambda))$.

*(a) If $u$ is a smonic morphism, then $v$ is a sepic morphism and $(*)$ is isomorphic to the following triangle:*

$$
\begin{array}{c}
X \\ \downarrow \\ Y \\ \downarrow \\ W \\ \downarrow \\ X[1]
\end{array}
\quad
\begin{array}{c}
\cdots \longrightarrow X^{i-1} \xrightarrow{d^{i-1}} X^i \xrightarrow{d^i} X^{i+1} \longrightarrow \cdots \\
\downarrow \binom{1}{0} \quad \downarrow \binom{1}{0} \quad \downarrow \binom{1}{0} \\
\cdots \longrightarrow X^{i-1} \oplus Y'^{(i-1)} \xrightarrow{\partial^{i-1}} X^i \oplus Y'^i \xrightarrow{\partial^i} X^{i+1} \oplus Y'^{(i+1)} \longrightarrow \cdots \\
\downarrow (0\ 1) \quad \downarrow (0\ 1) \quad \downarrow (0\ 1) \\
\cdots \longrightarrow Y'^{(i-1)} \xrightarrow{e^{i-1}} Y'^i \xrightarrow{e^i} Y'^{(i+1)} \longrightarrow \cdots \\
\downarrow -a^{i-1} \quad \downarrow -a^i \quad \downarrow -a^{i+1} \\
\cdots \longrightarrow X^i \xrightarrow{-d^i} X^{i+1} \xrightarrow{-d^{i+1}} X^{i+2} \longrightarrow \cdots
\end{array}
$$

*(b) If $u$ is a sepic morphism, then $v$ is a sirreducible morphism and $(*)$ is isomorphic to the following triangle:*

$$
\begin{array}{c}
X \\ \downarrow \\ Y \\ \downarrow \\ W \\ \downarrow \\ X[1]
\end{array}
\quad
\begin{array}{c}
\cdots \longrightarrow Y^{i-1} \xrightarrow{\partial^{i-1}} Y^i \xrightarrow{b^i} X'^{i+1} \xrightarrow{e^{i+1}} \cdots \\
\downarrow 1 \quad \downarrow 1 \quad \downarrow \\
\cdots \longrightarrow Y^{(i-1)} \xrightarrow{\partial^{i-1}} Y^i \longrightarrow 0 \longrightarrow \cdots \\
\quad \downarrow \quad \downarrow b^i \quad \downarrow \\
\cdots \longrightarrow 0 \longrightarrow X'^{(i+1)} \xrightarrow{-e^{i+1}} X'^{(i+2)} \longrightarrow \cdots \\
\quad \downarrow \quad \downarrow 1 \quad \downarrow 1 \\
\cdots \longrightarrow Y^i \xrightarrow{-b^i} X'^{i+1} \xrightarrow{-e^{i+1}} X'^{i+2} \longrightarrow \cdots
\end{array}
$$

*where $b^i$ is an irreducible morphism in $\mathcal{P}(\Lambda)$.*

*(c) If $u$ is a sirreducible morphism then:*

*(1) $v$ is a smonic morphism and $(*)$ is isomorphic to the following triangle:*

$$
\begin{array}{c}
X \\ \downarrow \\ Y \\ \downarrow \\ W \\ \downarrow \\ X[1]
\end{array}
\quad
\begin{array}{c}
\cdots \xrightarrow{e^{i-2}} X'^{(i-1)} \xrightarrow{e^{i-1}} X^i \longrightarrow 0 \longrightarrow \cdots \\
\downarrow \quad \downarrow \quad \downarrow f^i \quad \downarrow \\
\cdots \longrightarrow 0 \longrightarrow Y^i \xrightarrow{\epsilon^i} Y'^{(i+1)} \xrightarrow{\epsilon^{i+1}} \cdots \\
\quad \downarrow \quad \downarrow 1 \quad \downarrow 1 \\
\cdots \xrightarrow{-e^{i-1}} X^i \xrightarrow{f^i} Y^i \xrightarrow{\epsilon^i} Y'^{(i+1)} \longrightarrow \cdots \\
\quad \downarrow 1 \quad \downarrow \quad \downarrow \\
\cdots \xrightarrow{-e^{i-1}} X^i \longrightarrow 0 \longrightarrow 0 \longrightarrow \cdots
\end{array}
$$

*where $f^i$ is irreducible.*

*(2) $v$ is a sirreducible morphism and $(*)$ is isomorphic the following triangle:*

$$
\begin{array}{c}
X \\ \downarrow \\ Y \\ \downarrow \\ W \\ \downarrow \\ X[1]
\end{array}
\quad
\begin{array}{c}
\cdots \longrightarrow Y^{j-1} \longrightarrow Y^j \xrightarrow{b^j} X'^{j+1} \xrightarrow{e^{j+1}} \cdots \longrightarrow X'^{(i-1)} \xrightarrow{e^{i-1}} X^i \longrightarrow 0 \longrightarrow \cdots \\
\downarrow 1 \quad \downarrow 1 \quad \downarrow \quad \downarrow \quad \downarrow f^i \quad \downarrow \\
\cdots \longrightarrow Y^{j-1} \longrightarrow Y^j \longrightarrow 0 \longrightarrow \cdots \longrightarrow 0 \longrightarrow Y^i \xrightarrow{\epsilon^i} Y'^{i+1} \longrightarrow \cdots \\
\quad \downarrow \quad \downarrow b^j \quad \downarrow \quad \downarrow \quad \downarrow 1 \quad \downarrow 1 \\
\cdots \longrightarrow 0 \longrightarrow X'^{j+1} \xrightarrow{-e^{j+1}} X'^{j+2} \longrightarrow \cdots \xrightarrow{-e^{i-1}} X^i \xrightarrow{f^i} Y^i \xrightarrow{\epsilon^i} Y'^{i+1} \longrightarrow \cdots \\
\quad \downarrow \quad \downarrow 1 \quad \downarrow 1 \quad \downarrow \quad \downarrow \quad \downarrow \\
\cdots \longrightarrow Y^j \xrightarrow{-b^j} X'^{j+1} \xrightarrow{-e^{j+1}} X'^{j+2} \longrightarrow \cdots \xrightarrow{-e^{i-1}} X^i \longrightarrow 0 \longrightarrow 0 \longrightarrow \cdots
\end{array}
$$

*with some $b^j$ irreducible.*

*(3) $v$ is a a sirreducible morphism and and $(*)$ is isomorphic the following triangle:*

$$\begin{array}{c}
X \\ \downarrow \\ Y \\ \downarrow \\ W \\ \downarrow \\ X[1]
\end{array}
\quad
\begin{array}{c}
\cdots \longrightarrow Y^{i-3} \xrightarrow{\partial^{i-3}} Y^{i-2} \xrightarrow{\partial^{i-2}} Y^{i-1} \xrightarrow{c^{i-1}} X^i \longrightarrow 0 \longrightarrow \cdots \\
\downarrow 1 \quad\quad \downarrow 1 \quad\quad \downarrow 1 \quad\quad \downarrow f^i \quad\quad \downarrow \\
\cdots \longrightarrow Y^{i-3} \xrightarrow{\partial^{i-3}} Y^{i-2} \xrightarrow{\partial^{i-2}} Y^{i-1} \xrightarrow{\partial^{i-1}} Y^i \xrightarrow{\epsilon^i} Y'^{i+1} \xrightarrow{\epsilon^{i+1}} \cdots \\
\downarrow \quad\quad \downarrow \quad\quad \downarrow c^{i-1} \quad\quad \downarrow 1 \quad\quad \downarrow 1 \\
\cdots \longrightarrow 0 \longrightarrow 0 \longrightarrow X^i \xrightarrow{f^i} Y^i \xrightarrow{\epsilon^i} Y'^{i+1} \xrightarrow{\epsilon^{i+1}} \cdots \\
\downarrow \quad\quad \downarrow \quad\quad \downarrow 1 \quad\quad \downarrow \quad\quad \downarrow \\
\cdots \xrightarrow{-\partial^{i-3}} Y^{i-2} \xrightarrow{-\partial^{i-2}} Y^{i-1} \xrightarrow{-c^{i-1}} X^i \longrightarrow 0 \longrightarrow 0 \longrightarrow \cdots
\end{array}$$

with $c^{i-1}$ an irreducible morphism in $\mathcal{P}(\Lambda)$.

*Proof.* The proof of $(a)$ follows from Proposition (4).

$(b)$ Following Proposition (5) and Theorem (1), we have that if $u$ is sepic, then $v$ is sirreducible. So, supposing that $b^i$ is an irreducible morphism in $\mathcal{P}$ for some $i$, then $b^j$ is a split epimorphism for $j < i$ and $b^j$ is split monomorphism for $j > i$. Then, we can say that $X'^j = 0$ for $j \leq i$ and $Y^j = 0$ for $j \geq i+1$. So, the triangle $(*)$ is isomorphic to the triangle described in $(b)$.

$(c)$ Now, supposing that $u$ is sirreducible, then from Proposition (3) we have that $(*)$ is isomorphic to the following triangle

$$\begin{array}{c}
X \\ \downarrow \\ Y \\ \downarrow \\ W \\ \downarrow \\ X[1]
\end{array}
\quad
\begin{array}{c}
\cdots \longrightarrow Y^{(i-1)} \oplus X'^{(i-1)} \longrightarrow X^i \xrightarrow{d^i} X^{i+1} \xrightarrow{d^{i+1}} \cdots \\
\downarrow (1\,0) \quad\quad \downarrow f^i \begin{pmatrix}\ell^i \\ \epsilon^i\end{pmatrix} \quad\quad \downarrow \begin{pmatrix}1\\0\end{pmatrix} \\
\cdots \xrightarrow{\partial^{i-2}} Y^{(i-1)} \xrightarrow{\partial^{i-1}} Y^i \xrightarrow{} X^{i+1} \oplus Y'^{(i+1)} \longrightarrow \cdots \\
\downarrow c^{i-1} \quad\quad \downarrow 1 \quad\quad \downarrow (0\,1) \\
\cdots \xrightarrow{-e^{i-1}} X^i \xrightarrow{f^i} Y^i \xrightarrow{\epsilon^i} Y'^{(i+1)} \longrightarrow \cdots \\
\downarrow 1 \quad\quad \downarrow \ell^i \quad\quad \downarrow a^{i+1} \\
\cdots \xrightarrow{-(b^{i-1}\,e^{i-1})} X^i \xrightarrow{-d^i} X^{i+1} \xrightarrow{-d^{i+1}} X^{i+2} \longrightarrow \cdots
\end{array}$$

From Theorem (1), $v$ is isomorphic to a smonic or sirreducible morphism. If $v$ is isomorphic to a smonic morphism, then from Lemma (2), $v$ is smonic, so $c^{i-1}$ is a split monomorphism and $b^j$ is split monomorphism for $j < i-1$. So $Y^j = 0$ for $j \leq i-1$. We have too that $(0\,1) : X^j \oplus Y'^j \to Y'^j$ is a monomorphism for $j \geq i+1$. So, $X^j = 0$ for $j \geq i+1$. Then, in this case, the previous triangle is isomorphic to the triangle described in $(c)(1)$.

Now, supposing that $v$ is isomorphic to a sirreducible morphism, we have that for some $j$, $b^j$ is an irreducible morphism for some $j \leq i-2$ or $c^{i-1}$ is an irreducible morphism. Supposing that $b^j$ is an irreducible morphism for some $j \leq i-2$.

$$\begin{array}{c}
X \\ \downarrow \\ Y \\ \downarrow \\ W \\ \downarrow \\ X[1]
\end{array}
\quad
\begin{array}{c}
\cdots \quad Y^{j-1} \oplus X'^{j-1} \longrightarrow Y^j \oplus X'^j \xrightarrow{\begin{pmatrix}\partial^j & 0 \\ b^j & e^j\end{pmatrix}} Y^{j+1} \oplus X'^{j+1} \longrightarrow \cdots \\
\downarrow (1\,0) \quad\quad \downarrow (1\,0) \quad\quad \downarrow (1\,0) \\
\cdots \quad Y^{j-1} \longrightarrow Y^j \xrightarrow{\partial^j} Y^{j+1} \longrightarrow \cdots \\
\downarrow b^{j-1} \quad\quad \downarrow b^j \quad\quad \downarrow b^{j+1} \\
\cdots \quad X'^j \longrightarrow X'^{j+1} \xrightarrow{-e^{j+1}} X'^{j+2} \longrightarrow \cdots \\
\downarrow \begin{pmatrix}1\\0\end{pmatrix} \quad\quad \downarrow \begin{pmatrix}1\\0\end{pmatrix} \quad\quad \downarrow \begin{pmatrix}1\\0\end{pmatrix} \\
\cdots \quad Y^j \oplus X'^j \longrightarrow Y^{j+1} \oplus X'^{j+1} \longrightarrow Y^{j+2} \oplus X'^{j+2} \longrightarrow \cdots
\end{array}$$

Then, $b^k : Y^k \to X'^{k+1}$ is a split epimorphism for $k < j$. So, $X'^{k+1} = 0$ for $k < j$. We have that $b^k$ is a split monomorphism for $k > j$. So, $Y^k = 0$ for $k > j$. Then, the previous diagram can be rewritten as follow:

$$\begin{array}{c}
X \\ \downarrow \\ Y \\ \downarrow \\ W \\ \downarrow \\ X[1]
\end{array}
\quad
\begin{array}{c}
\cdots \quad Y^{j-1} \longrightarrow Y^j \xrightarrow{b^j} X'^{j+1} \xrightarrow{e^{j+1}} \cdots \\
\downarrow 1 \quad\quad \downarrow 1 \quad\quad \downarrow \\
\cdots \quad Y^{j-1} \longrightarrow Y^j \longrightarrow 0 \longrightarrow \cdots \\
\downarrow \quad\quad \downarrow b^j \quad\quad \downarrow \\
\cdots \quad 0 \longrightarrow X'^{j+1} \xrightarrow{-e^{j+1}} X'^{j+2} \longrightarrow \cdots \\
\downarrow \quad\quad \downarrow 1 \quad\quad \downarrow 1 \\
\cdots \quad Y^j \xrightarrow{-b^j} X'^{j+1} \xrightarrow{-e^{j+1}} X'^{j+2} \longrightarrow \cdots
\end{array}$$

We also have that $c^{i-1} : Y^{i-1} \to X^i$ is a split monomorphism and $(0\,1) : X^j \oplus Y'^j \to Y'^j$ is a split monomorphism for $j \geq i+1$ in the following diagram:

$$\begin{array}{c}
X \\
\downarrow \\
Y \\
\downarrow \\
W \\
\downarrow \\
X[1]
\end{array}
\quad
\begin{array}{ccccccccc}
\cdots & \longrightarrow & Y^{(i-1)} \oplus X'^{(i-1)} & \longrightarrow & X^i & \xrightarrow{d^i} & X^{i+1} & \xrightarrow{d^{i+1}} & \cdots \\
& & \downarrow{(1\ 0)} & & \downarrow f^i \binom{\ell^i}{\epsilon^i} & & \downarrow \binom{1}{0} & & \\
\cdots & \xrightarrow{\partial^{i-2}} & Y^{(i-1)} & \xrightarrow{\partial^{i-1}} & Y^i & \xrightarrow{} & X^{i+1} \oplus Y'^{(i+1)} & \longrightarrow & \cdots \\
& & \downarrow c^{i-1} & & \downarrow 1 & & \downarrow (0\ 1) & & \\
\cdots & \xrightarrow{-e^{i-1}} & X^i & \xrightarrow{f^i} & Y^i & \xrightarrow{\epsilon^i} & Y'^{(i+1)} & \longrightarrow & \cdots \\
& & \downarrow 1 & & \downarrow \ell^i & & \downarrow a^{i+1} & & \\
\cdots & \xrightarrow{-(b^{i-1}\ e^{i-1})} & X^i & \xrightarrow{-d^i} & X^{i+1} & \xrightarrow{-d^{i+1}} & X^{i+2} & \longrightarrow & \cdots
\end{array}$$

So, $Y^{i-1} = 0$ and $X^j = 0$ for $j \geq i + 1$. So, we can simplify the diagram and obtain the following diagram:

$$\begin{array}{c}
X \\
\downarrow \\
Y \\
\downarrow \\
W \\
\downarrow \\
X[1]
\end{array}
\quad
\begin{array}{ccccccccc}
\cdots & \xrightarrow{e^{i-2}} & X'^{(i-1)} & \xrightarrow{e^{i-1}} & X^i & \longrightarrow & 0 & \longrightarrow & \cdots \\
& & \downarrow & & \downarrow f^i & & \downarrow & & \\
\cdots & \longrightarrow & 0 & \longrightarrow & Y^i & \xrightarrow{\epsilon^i} & Y'^{(i+1)} & \xrightarrow{\epsilon^{i+1}} & \cdots \\
& & \downarrow & & \downarrow 1 & & \downarrow 1 & & \\
\cdots & \xrightarrow{-e^{i-1}} & X^i & \xrightarrow{f^i} & Y^i & \xrightarrow{\epsilon^i} & Y'^{(i+1)} & \longrightarrow & \cdots \\
& & \downarrow 1 & & \downarrow & & \downarrow & & \\
\cdots & \xrightarrow{-e^{i-1}} & X^i & \longrightarrow & 0 & \longrightarrow & 0 & \longrightarrow & \cdots
\end{array}$$

Then $(*)$ is isomorphic to the triangle described in $(c)(2)$.

Now, supposing that $c^{i-1}$ is irreducible. Then, $b^j$ is split epimorphism for $j \leq i-2$. So $X'^j = 0$ for $j \leq i-1$. And we also have that $(0\ 1) : X^j \oplus Y'^j \to Y'^j$ is a split monomorphism for each $j \geq i+1$. So, $X^j = 0$ for $j \geq i+1$. Then, the Auslander-Reiten triangle $(*)$ is isomorphic to the triangle describe in $(c)(3)$.

$\square$

## 5. The mapping cone of an irreducible morphism

One of the objectives of this article is to characterize the mapping cone of an Auslander-Reiten triangle. It is possible to prove that this object is indecomposable. We are able to prove it in a more general context, in the triangulated category. In [16], Propositon 6.1, proved that the mapping cone of irreducible morphism between indecomposable object is indecomposable in $D^b(\text{mod}\ \Lambda)$. So, irreducible morphisms provide a useful method to construct indecomposable objects in a derived category. It is important to note that, not all irreducible morphisms appear in Auslander-Reiten triangles (see [24]).

### 5.1. The mapping cone of an irreducible morphism is indecomposable.
Now we are supposing that $\mathcal{T}$ is a triangulated category.

**Proposition 6.** *Let $X$ and $Y$ be indecomposable in a triangulated category $\mathcal{T}$ and assume that we have an irreducible morphism $u : X \to Y$. Then, the mapping cone $C_u$ is indecomposable.*

*Proof.* Let $X \xrightarrow{u} Y \xrightarrow{v} C_u \xrightarrow{w} X[1]$ be a distinguished triangle in $\mathcal{T}$ and $f : C_u \to C_u$ an idempotent. We have to prove that $f$ is zero or the identity morphism to conclude that $C_u$ is indecomposable. Let $C_u \xrightarrow{f} C_u \xrightarrow{g} C_f \xrightarrow{h} C_u[1]$ be a triangle and then, from the octaedral axiom, we have the following commutative diagram

$$\begin{array}{ccccccc}
C_u & \xrightarrow{f} & C_u & \xrightarrow{g} & C_f & \xrightarrow{h} & C_u[1] \\
\downarrow 1 & & \downarrow w & & \downarrow \alpha & & \downarrow 1 \\
C_u & \xrightarrow{wf} & X[1] & \xrightarrow{\theta} & U & \xrightarrow{\delta} & C_u[1] \\
& & \downarrow -u[1] & & \downarrow \beta & & \downarrow f[1] \\
& & Y[1] & \xrightarrow{1} & Y[1] & \xrightarrow{-v[1]} & C_u[1] \\
& & \downarrow -v[1] & & \downarrow \gamma & & \\
& & C_u[1] & \xrightarrow{g[1]} & C_f[1] & &
\end{array}$$

So, from $-u[1] = \beta\theta$ and $u$ irreducible, we have that $\theta$ is a section or $\beta$ is a retraction.

If $\theta$ is a section, then $wf = 0$, so there is $\eta : C_u \to Y$ such that $v\eta = f$. From $Y$ indecomposable, then $\eta v : Y \to Y$ is an isomorphism or nilpotent. If $\eta v$ is an isomorphism, then $v$ is a section and then $u = 0$ and this is a contradiction. If $\eta v$ is nilpotent such that $(\eta v)^n = 0$ and $(\eta v)^{n-1} \neq 0$, then from $f^n = f$ we have $f = 0$. Similarly, if $\beta$ is a retraction and then $u = 0$ and we have a contradiction. □

5.2. **The mapping cone of an irreducible morphism that starts in an object of** $\mod \Lambda$**.** Here we present a general result about irreducible morphism in $D^b(\mod \Lambda)$ that starts in an object $X \in \mod \Lambda$.

**Proposition 7.** *Let* $f : X \to Y$ *be an irreducible morphism in* $D^b(\mod \Lambda)$

*(a) If* $X \in \mod A$*, then* $\mathrm{Hom}_{D^b(\mod \Lambda)}(C_f, P) = 0$ *for all indecomposable projective* $\Lambda$*-module.*
*(b) If* $Y \in \mod \Lambda$*, then* $\mathrm{Hom}_{D^b(\mod \Lambda)}(I, C_f[-1]) = 0$ *for all indecomposable injective* $\Lambda$*-module.*

*Proof.* $(a)$ Let $X \xrightarrow{f} Y \xrightarrow{g} C_f \xrightarrow{h} X[1]$ be the triangle associated with $f$ and let $u : C_f \to P$ be a non zero morphism in $D^b(\mod A)$. Since $Ext_A^{-1}(X, P) = 0$, then $ug \neq 0$. From the commutative diagram

$$\begin{array}{ccccccc} Y & \xrightarrow{g} & C_f & \xrightarrow{h} & X[1] & \xrightarrow{-f[1]} & Y[1] \\ \downarrow 1 & & \downarrow u & & \downarrow \theta & & \downarrow 1 \\ Y & \xrightarrow{ug} & P & \xrightarrow{w} & C_{ug} & \xrightarrow{\alpha} & Y[1] \end{array}$$

and $f$ irreducible implies either $\theta$ is a split monomorphism or $\alpha$ is a split epimorphism. Since $ug \neq 0$ if and only if $\alpha$ is not a split epimorphism, then $\theta$ is a split monomorphism. So, from $\theta h = wu$, we have $h = \theta' wu$, where $\theta'\theta = 1$. So $\theta' w \in \mathrm{Hom}_{D^b(\mod \Lambda)}(P, X[1]) = Ext_A^1(P, X) = 0$ and then $h = 0$. This implies that $f$ is a split monomorphism.

$(b)$ The proof follows in a similar way.
□

5.3. **Conditions to the non existence of irreducible morphisms.**

**Proposition 8.** *Let* $X, Y$ *be complexes such that,* $X^i = Y^i = 0$ *for all* $i > 0$ *and there exists a* $t > 0$ *such that* $X^{-i} = 0$ *for all* $i > t$ *and* $X^{-t} \neq 0$*. If* $f : X \to Y$ *is a irreducible morphism, then* $Y^{-k} = 0$ *for all* $k \geq t + 2$*.*

*Proof.* Suppose that $Y^{-k} \neq 0$ for some $k \geq t+2$. We have the following factorisation of $f = gh$.

$$\begin{array}{ccccccccccc} \cdots & \to & 0 & \to & 0 & \to & X^{-t} & \to & X^{-t+1} & \to & \cdots & \to & X^0 & \to & 0 & \to \cdots \\ & & \downarrow & & \downarrow & & \downarrow f^{-t} & & \downarrow f^{-t+1} & & & & \downarrow f^0 & & \downarrow & \\ \cdots & \to & 0 & \to & Y^{(-t-1)} & \to & Y^{-t} & \to & Y^{-t+1} & \to & \cdots & \to & Y^0 & \to & 0 & \to \cdots \\ & & \downarrow & & \downarrow 1 & & \downarrow 1 & & \downarrow 1 & & & & \downarrow 1 & & \downarrow & \\ \cdots & \to & Y^{(-t-2)} & \to & Y^{(-t-1)} & \to & Y^{-t} & \to & Y^{-t+1} & \to & \cdots & \to & Y^0 & \to & 0 & \to \cdots \end{array}$$

Then, $h$ is not a split monomorphism, since $f$ is not a split monomorphism and since $Y^k \neq 0$ then $g$ is not a split epimorphism. So, $f$ is not irreducible. □

So, we have the following result of Happel-Keller-Reiten (Proposition 6.2 - [16]).

**Corollary 1.** *Let* $\Lambda$ *be a finite dimensional algebra,* $X$ *and* $Y$ *indecomposable* $\Lambda$*-modules with* $\mathrm{pd}_\Lambda X < \infty$ *and* $\mathrm{pd}_\Lambda Y \geq \mathrm{pd}_\Lambda X + 2$. *Then, there is no irreducible morphism* $f : X \to Y$ *in* $D^b(\mod \Lambda)$*.*

6. A CONSTRUCTION OF A COMPONENT OF THE AUSLANDER-REITEN QUIVER OF THE DERIVED CATEGORY

In [1], Assem and Brenner, had developed an algorithm in which it is possible to construct a tower of full additive subcategories of the homotopy category $\mathcal{K}^b(\mathcal{P}(\Lambda))$ of bounded complexes of projective objects in $\mod \Lambda$. It offers some help to construct the Aulander-Reiten quiver of the derived category of some tilted algebras. In this section, we are going to describe, using two examples, an explicit construction of the component of the Aulander-Reiten quiver of the derived category of a tilted algebra that contains the slice. This technique can be applied in general to the construction of every component of the Auslander-Reiten quiver of the derived category that contains the slice.

We know that irreducible morphisms in module category do not generally stay irreducible in the derived category of the module category. Now, we describe a specific situation of irreducible morphisms in the module category that can remain irreducible in the derived category and this will help us construct the component of the Auslander-Reiten quiver of the derived category that contains this specific irreducible morphisms.

Besides, after the construction of the component, it will be easy to visualize where the module category $\mod \Lambda$ is. In order to do this, it is necessary to know the representation of the projective resolution of each module of the slice. We refer to [17] and [22] for basic facts about tilted algebras.

Let $\Lambda$ be a finite dimensional hereditay algebra over a field $k$ and $\mod \Lambda$ be the category of finitely generated left $\Lambda$-modules. We denote by $D = \Hom_k(-, k)$ the duality between left and right $A$-modules. Given a tilting module $T \in \mod \Lambda$ with $\Gamma = \End_\Lambda(T)$, we call $\Gamma$ a tilted algebra and we denote by $\mathcal{T}$ the full subcategory of $\mod \Lambda$ whose objects are generated by $T$ and by $\mathcal{F}$ the full subcategory of $\mod \Lambda$ whose objects $X$ satisfy $\Hom_\Lambda(T, X) = 0$. It is possible to prove that $\mathcal{T}$ is the full subcategory of $\mod \Lambda$ whose objects $X$ satisfy $\Ext^1_\Lambda(T, X) = 0$ and that $(\mathcal{T}, \mathcal{F})$ is a torsion pair in $\mod \Lambda$.

Let $D(\Lambda_\Lambda)$ be the injective cogenerator of $\mod \Lambda$. It is possible to prove that $D(\Lambda_\Lambda)$ is contained in the subcategory $\mathcal{T}$ and that $\Hom_\Lambda(T, D(\Lambda_\Lambda))$ is a complete slice in $\mod \Gamma$. Accordingly to Happel ([13]), the derived functor $R\Hom_\Lambda(T, -) : D^b(\mod \Lambda) \to D^b(\mod \Gamma)$ gives us a derived equivalence such that $(\mathcal{X}, \mathcal{Y})$ is a splitting torsion pair in $\mod \Gamma$ with $\mathcal{Y} = F(\mathcal{T})$ and $\mathcal{X} = F(\mathcal{F}[1])$ and we have $R\Hom_\Lambda(T, D(\Lambda_\Lambda)) = \Hom_\Lambda(T, D(\Lambda_\Lambda))$. So, for each irreducible morphism $I \to I'$ in $\mod \Lambda$ with $I, I'$ indecomposable injective modules, we have the corresponding irreducible morphim $I \to I'$ in $D^b(\mod \Lambda)$ and so we have an irreducible morphism $R\Hom_\Lambda(T, I) \to R\Hom_\Lambda(T, I')$ in $D^b(\mod \Gamma)$.

(1) Let $\Lambda$ be the hereditary path algebra of the quiver $Q$

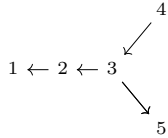

of type $\mathbb{D}_5$. Let $T = P_5 \oplus P_4 \oplus \tau^{-2} P_5 \oplus I_5 \oplus I_4$, where $P_4, P_5$ are the indecomposable projective modules and $I_4, I_5$ are the indecomposable injective modules associated to the vertices $4$ and $5$. It is easily checked that $T$ is a tilting $\Lambda$-module and $\Gamma = \End_\Lambda(T)$ is a tilted algebra given by the quiver

$$1 \xleftarrow{\delta} 2 \xleftarrow{\gamma} 3 \xleftarrow{\beta} 4 \xleftarrow{\alpha} 5$$

bounded by $\alpha\beta\gamma\delta = 0$. Computing the Auslander-Reiten quiver of $\Gamma$ yields

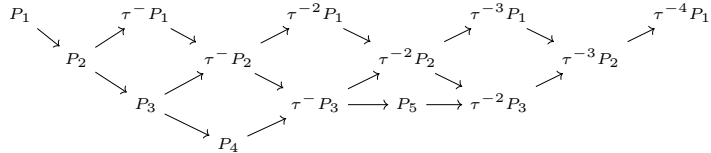

We have that $\tau^{-1} P_1, \tau^- P_2, \tau^- P_3, P_5$ and $P_4$ is a complete slice in $\mod \Gamma$ with the following minimal projective resolutions: $0 \to P_1 \to P_2 \to \tau^{-1} P_1 \to 0$, $0 \to P_1 \to P_3 \to \tau^{-1} P_2 \to 0$, $0 \to P_1 \to P_4 \to \tau^{-1} P_3 \to 0$. So, in the homotopic category $\mathcal{K}^b(\mathcal{P}(\Gamma))$ we represent each one of these projective resolutions by $(P_1 - P_2[0]), (P_1 - P_3[0]), (P_1 - P_4[0]), P_4[0], P_5[0]$. So, the strategy for constructing the component is: taking each morphism in the slice, construct the mapping cone of each one and so on applying the Theorem (2).

We have the following component of the Auslander-Reiten quiver of the bounded derived category. In the picture below we have only represented the objects that represent a module of $\mod \Gamma$ (represented by its projective resolutions).

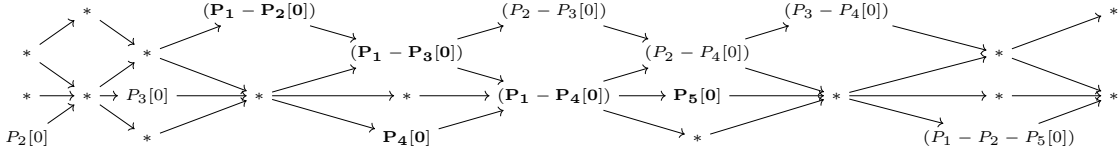

(2) Let $\Lambda$ be the representation infinite hereditary path algebra of the quiver $Q$

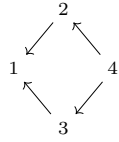

and $T = P_1 \oplus P_4/P_2 \oplus P_4/P_3 \oplus I_4$. $T$ is a tilting $\Lambda$-module and $\Gamma = \mathrm{End}_\Lambda(T)$ is a tilted algebra given by the quiver

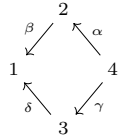

with relations $\alpha\beta = 0 = \gamma\delta$. By computing the Auslander-Reiten quiver of $\Gamma$, yelds

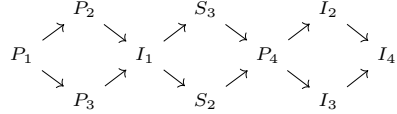

The objects $I_1, S_2, S_3$ and $P_4$ form a complete slice in $\mathrm{mod}\,\Gamma$ and each one will be represented by its respective projective resolution in the following quiver

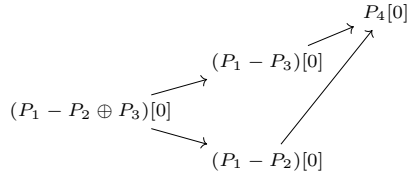

So, in doing each respective cone and applying the Theorem (2) we have the transjective component of the Auslander-Reiten quiver of the derived category of $\mathrm{mod}\,\Gamma$

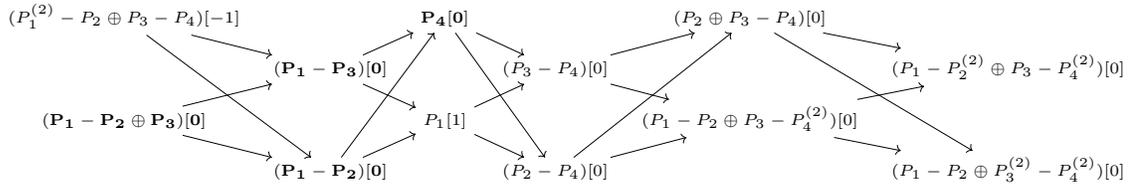

In this way it is possible to note that the only modules of $\mathrm{mod}\,\Gamma$ that appear in this transjective component are those modules in the slice.

## Acknowledgement

The authors wish to thank Patrick Le Meur and Ibrahim Assem for their fruitful discussions. The first author was partially supported by the DMAT-UFPR and CNPq-Universal 477880/2012-6.The second author was partially supported by DMA-UFV. The third author was partially supported by CODI and Estrategia de Sostenibilidad 2016-2017 (Universidad de Antioquia), and COLCIENCIAS-ECOPETROL (Contrato RC. No. 0266-2013).

(Edson Ribeiro Alvares) CENTRO POLITÉCNICO, DEPARTAMENTO DE MATEMÁTICA, UNIVERSIDADE FEDERAL DO PARANÁ, CP019081, JARDIM DAS AMERICAS, CURITIBA-PR, 81531-990, BRAZIL
*E-mail address*: rolo1rolo@gmail.com, rolo@ufpr.br

(Sônia Maria Fernandes) DEPARTAMENTO DE MATEMÁTICA, UNIVERSIDADE FEDERAL DE VIÇOSA, VIÇOSA, BRAZIL.
*E-mail address*: somari@ufv.br

(Hernán Giraldo) INSTITUTO DE MATEMÁTICAS, UNIVERSIDAD DE ANTIOQUIA, CALLE 67 NO. 53-108, MEDELLÍN, COLOMBIA.
*E-mail address*: hernan.giraldo@udea.edu.co